\documentclass[10pt,a4paper]{article}
\usepackage{amsfonts}
\usepackage{amssymb}
\usepackage{mathrsfs}
\usepackage{amsthm}
\usepackage{setspace}
\usepackage{color}

\newtheorem{theo}{Theorem}[section]
\newtheorem{prop}[theo]{Proposition}

\newcommand{\be}{\begin{eqnarray*}}
\newcommand{\ee}{\end{eqnarray*}}
\newcommand{\beqa}{\begin{eqnarray}}
\newcommand{\eeqa}{\end{eqnarray}}
\newcommand{\ba}{\begin{array}}
\newcommand{\ea}{\end{array}}

\newcommand{\mf}{\mathfrak}
\newcommand{\mbb}{\mathbb}

\newcommand{\BasicTri}{\begin{picture}(226,156)%
\put(4,13){\line(1,0){216}}%
\put(4,141){\line(1,0){216}}%
\put(36,77){\line(1,0){152}}%
\put(10,1){\line(1,2){76}}%
\put(74,1){\line(1,2){76}}%
\put(138,1){\line(1,2){76}}%
\put(86,1){\line(-1,2){76}}%
\put(150,1){\line(-1,2){76}}%
\put(214,1){\line(-1,2){76}}%
\put(3,7){\makebox(0,0){\phantom{p}}}%
\end{picture}}

\newcommand{\TriAll}{\begin{picture}(226,156)%
\put(4,13){\line(1,0){216}}%
\put(4,141){\line(1,0){216}}%
\put(36,77){\line(1,0){152}}%
\put(10,1){\line(1,2){76}}%
\put(74,1){\line(1,2){76}}%
\put(138,1){\line(1,2){76}}%
\put(86,1){\line(-1,2){76}}%
\put(150,1){\line(-1,2){76}}%
\put(214,1){\line(-1,2){76}}%

\put(47,141){\line(-1,-2){5}}\put(49,141){\line(-1,-2){5}}\put(47,141){\line(-1,2){5}}\put(49,141){\line(-1,2){5}}%
\put(47,13){\line(-1,-2){5}}\put(49,13){\line(-1,-2){5}}\put(47,13){\line(-1,2){5}}\put(49,13){\line(-1,2){5}}%
\put(175,141){\line(-1,-2){5}}\put(177,141){\line(-1,-2){5}}\put(175,141){\line(-1,2){5}}\put(177,141){\line(-1,2){5}}%
\put(175,13){\line(-1,-2){5}}\put(177,13){\line(-1,-2){5}}\put(175,13){\line(-1,2){5}}\put(177,13){\line(-1,2){5}}%
\put(143,77){\line(-1,-2){5}}\put(145,77){\line(-1,-2){5}}\put(143,77){\line(-1,2){5}}\put(145,77){\line(-1,2){5}}%

\put(111,13){\line(1,-2){5}}\put(113,13){\line(1,-2){5}}\put(111,13){\line(1,2){5}}\put(113,13){\line(1,2){5}}%
\put(111,141){\line(1,-2){5}}\put(113,141){\line(1,-2){5}}\put(111,141){\line(1,2){5}}\put(113,141){\line(1,2){5}}%
\put(79,77){\line(1,-2){5}}\put(81,77){\line(1,-2){5}}\put(79,77){\line(1,2){5}}\put(81,77){\line(1,2){5}}%

\put(31,43){\line(-1,0){10}}\put(31,43){\line(1,-2){5}}\put(32,45){\line(-1,0){10}}\put(32,45){\line(1,-2){5}}%
\put(95,43){\line(-1,0){10}}\put(95,43){\line(1,-2){5}}\put(96,45){\line(-1,0){10}}\put(96,45){\line(1,-2){5}}%
\put(159,43){\line(-1,0){10}}\put(159,43){\line(1,-2){5}}\put(160,45){\line(-1,0){10}}\put(160,45){\line(1,-2){5}}%
\put(63,107){\line(-1,0){10}}\put(63,107){\line(1,-2){5}}\put(64,109){\line(-1,0){10}}\put(64,109){\line(1,-2){5}}%
\put(127,107){\line(-1,0){10}}\put(127,107){\line(1,-2){5}}\put(128,109){\line(-1,0){10}}\put(128,109){\line(1,-2){5}}%
\put(191,107){\line(-1,0){10}}\put(191,107){\line(1,-2){5}}\put(192,109){\line(-1,0){10}}\put(192,109){\line(1,-2){5}}%

\put(64,46){\line(-1,-2){5}}\put(64,46){\line(1,0){10}}\put(65,44){\line(-1,-2){5}}\put(65,44){\line(1,0){10}}%
\put(128,46){\line(-1,-2){5}}\put(128,46){\line(1,0){10}}\put(129,44){\line(-1,-2){5}}\put(129,44){\line(1,0){10}}%
\put(192,46){\line(-1,-2){5}}\put(192,46){\line(1,0){10}}\put(193,44){\line(-1,-2){5}}\put(193,44){\line(1,0){10}}%
\put(32,110){\line(-1,-2){5}}\put(32,110){\line(1,0){10}}\put(33,108){\line(-1,-2){5}}\put(33,108){\line(1,0){10}}%
\put(96,110){\line(-1,-2){5}}\put(96,110){\line(1,0){10}}\put(97,108){\line(-1,-2){5}}\put(97,108){\line(1,0){10}}%
\put(160,110){\line(-1,-2){5}}\put(160,110){\line(1,0){10}}\put(161,108){\line(-1,-2){5}}\put(161,108){\line(1,0){10}}%

\end{picture}}

\newcommand{\TriAcross}{\begin{picture}(226,156)%
\put(4,13){\line(1,0){216}}%
\put(4,141){\line(1,0){216}}%
\put(36,77){\line(1,0){152}}%
\put(10,1){\line(1,2){76}}%
\put(74,1){\line(1,2){76}}%
\put(138,1){\line(1,2){76}}%
\put(86,1){\line(-1,2){76}}%
\put(150,1){\line(-1,2){76}}%
\put(214,1){\line(-1,2){76}}%

\put(47,141){\line(-1,-2){5}}\put(49,141){\line(-1,-2){5}}\put(47,141){\line(-1,2){5}}\put(49,141){\line(-1,2){5}}%
\put(47,13){\line(-1,-2){5}}\put(49,13){\line(-1,-2){5}}\put(47,13){\line(-1,2){5}}\put(49,13){\line(-1,2){5}}%
\put(175,141){\line(-1,-2){5}}\put(177,141){\line(-1,-2){5}}\put(175,141){\line(-1,2){5}}\put(177,141){\line(-1,2){5}}%
\put(175,13){\line(-1,-2){5}}\put(177,13){\line(-1,-2){5}}\put(175,13){\line(-1,2){5}}\put(177,13){\line(-1,2){5}}%
\put(143,77){\line(-1,-2){5}}\put(145,77){\line(-1,-2){5}}\put(143,77){\line(-1,2){5}}\put(145,77){\line(-1,2){5}}%

\put(111,13){\line(1,-2){5}}\put(113,13){\line(1,-2){5}}\put(111,13){\line(1,2){5}}\put(113,13){\line(1,2){5}}%
\put(111,141){\line(1,-2){5}}\put(113,141){\line(1,-2){5}}\put(111,141){\line(1,2){5}}\put(113,141){\line(1,2){5}}%
\put(79,77){\line(1,-2){5}}\put(81,77){\line(1,-2){5}}\put(79,77){\line(1,2){5}}\put(81,77){\line(1,2){5}}%

\end{picture}}

\newcommand{\TriUp}{\begin{picture}(226,156)%
\put(4,13){\line(1,0){216}}%
\put(4,141){\line(1,0){216}}%
\put(36,77){\line(1,0){152}}%
\put(10,1){\line(1,2){76}}%
\put(74,1){\line(1,2){76}}%
\put(138,1){\line(1,2){76}}%
\put(86,1){\line(-1,2){76}}%
\put(150,1){\line(-1,2){76}}%
\put(214,1){\line(-1,2){76}}%

\put(31,43){\line(-1,0){10}}\put(31,43){\line(1,-2){5}}\put(32,45){\line(-1,0){10}}\put(32,45){\line(1,-2){5}}%
\put(95,43){\line(-1,0){10}}\put(95,43){\line(1,-2){5}}\put(96,45){\line(-1,0){10}}\put(96,45){\line(1,-2){5}}%
\put(159,43){\line(-1,0){10}}\put(159,43){\line(1,-2){5}}\put(160,45){\line(-1,0){10}}\put(160,45){\line(1,-2){5}}%
\put(63,107){\line(-1,0){10}}\put(63,107){\line(1,-2){5}}\put(64,109){\line(-1,0){10}}\put(64,109){\line(1,-2){5}}%
\put(127,107){\line(-1,0){10}}\put(127,107){\line(1,-2){5}}\put(128,109){\line(-1,0){10}}\put(128,109){\line(1,-2){5}}%
\put(191,107){\line(-1,0){10}}\put(191,107){\line(1,-2){5}}\put(192,109){\line(-1,0){10}}\put(192,109){\line(1,-2){5}}%

\put(64,46){\line(-1,-2){5}}\put(64,46){\line(1,0){10}}\put(65,44){\line(-1,-2){5}}\put(65,44){\line(1,0){10}}%
\put(128,46){\line(-1,-2){5}}\put(128,46){\line(1,0){10}}\put(129,44){\line(-1,-2){5}}\put(129,44){\line(1,0){10}}%
\put(192,46){\line(-1,-2){5}}\put(192,46){\line(1,0){10}}\put(193,44){\line(-1,-2){5}}\put(193,44){\line(1,0){10}}%
\put(32,110){\line(-1,-2){5}}\put(32,110){\line(1,0){10}}\put(33,108){\line(-1,-2){5}}\put(33,108){\line(1,0){10}}%
\put(96,110){\line(-1,-2){5}}\put(96,110){\line(1,0){10}}\put(97,108){\line(-1,-2){5}}\put(97,108){\line(1,0){10}}%
\put(160,110){\line(-1,-2){5}}\put(160,110){\line(1,0){10}}\put(161,108){\line(-1,-2){5}}\put(161,108){\line(1,0){10}}%

\end{picture}}

\begin{document}

\title{A note on discrete Holonomy through directed edges, with no lengths}
\author{Stuart Armstrong and Jussi Westergren}
\date{2009}
\maketitle

\begin{abstract}
This note demonstrates how both the concept of distance and the concept of holonomy can be constructed from a suitable network with directed edges (and no lengths). The number of different edge types depends on the signature of the metric and the dimension of the holonomy group. If the holonomy group is of dimension one and the metric is positive-definite, a single type of directed edges is needed.
\end{abstract}

\section{Introduction}

Replacing continuous manifolds with discreet simplicial complexes is a very fertile area of current research - consider Regge calculus \cite{regge}, simplicial manifolds, and Seth Lloyds work on quantum gravity on a lattice derived from quantum computations \cite{sethGrav}. An extension of these current results is the use of causal sets (see Joe Henson's paper \cite{hensonGrav} for an overview and an application to quantum gravity), and the construction of gauge theories on them dependent on certain holonomy groups (\cite{romanBoson}, \cite{romanSpinor} and \cite{romanGauge}).

The causal set constructions are very interesting, as they deduce the whole structure -- the metric, time- and space-like distances, and so on -- from only the causal links, with no extra information. The other constructions have a large amount of extraneous information, though -- the lengths of the edges, or the holonomy group itself, which is a continuous valued transformation and very much against the spirit of discrete mathematics.

This paper will demonstrate how all these constructions can be replaced with a simple network combined with a small collection of different type of directed edges. If only positive distances can exist in the network, we shall need $d$ types of directed edges, where $d$ is the dimension of the holonomy group. If there are negative distances allowed in the network (for instance, a simplicial complex based on an underlying split-signature metric), we shall need $2d$ types of directed edges (or, equivalently, $d$ directed edges and $2$ undirected edges,

The construction is not gauge invariant.

\section{Length from edges}
It is very easy to replace a simplicial complex with edges-with-lengths with another complex with simple edges -- if there is a scale $\epsilon$ below which variations in lengths have no observable large scale consequences. Simply replace each edge of length $l$ with $m$ lengthless edges, where $m$ is
\be
m = \lfloor l / (2 \times \epsilon) \rfloor.
\ee
The requirement that length changes below $\epsilon$ have no observable consequences is what one would expect in a quantum universe. Then by seeing each edge as having the same length $\epsilon / 2$, the whole distance functions, geodesic lengths, etc... can be reconstructed.

If some of the lengths are negative, we shall need to use two types of edges: one for positive distances, and one for the negative distances.

\section{Holonomy from edges: $U(1)$ holonomy}
There are alternative ways of getting holonomy on a discrete set, if we are prepared to give up invariance. For the purpose of illustration, let $E \to M$ be a complex line bundle with connection $\nabla$, with holonomy contained in $U(1)$-- in other words, locally $E \cong \mbb{C} \times M$, and parallel transport via $\nabla$ can only change the phase (i.e. the difference between two parallel transports is given by a unit complex number).

\begin{prop}
If we choose a never zero section $s$ of $E$, we have a canonical identification $E = \mbb{C} \times M$.
\end{prop}
\begin{proof}
Let $x \in E$ be given. Then the fibre $E_x$ of $E$ at $x$ is a one dimensional complex space, and $s(x)$ is a non-zero element of it. Therefore, for any point $t \in E_x$, $t = \lambda s(x)$ for some unique complex number $\lambda$. Hence we may identify $t$ with $(x,\lambda)$; varying $t$ and $x$ then gives the canonical identification $E = \mbb{C} \times M$.
\end{proof}

This choice of $s$ (it can also be a local choice of $s$, if there are no global sections) is the invariance we give up in this case. The choice of $s$ is arbitrary, but, once it is chosen, the construction proceeds naturally. The main point is that we now have a canonical concept of the phase change along a non-closed path: simply parallel transport $s$ and compute the phase difference from the original $s$ at the end-point.

\subsection{Phase changes from directed edges}

The idea is as follows. First, cut up the space into simplexes, and leave the edges undirected. Then the distance function on this space is simply given by the construction above.

Then direct the edges. Parallel transport along an edge is then interpreted as being a change of phase by a certain (fixed) amount. Reversing the direction is interpreted as being a change of phase by a the inverse amount. For this to work, we need there to be an $\epsilon'>0$ such that the phase change $e^{i\epsilon'}$ has no observable consequences.

By changing the direction of the edges, we can construct all types of holonomy, though we may need to go to a finer subdivision than the ones previously.

For example, imagine we start with the plane, subdivided into triangles (this is not particularly accurate for geodesic distance, but is very useful for picturing things).
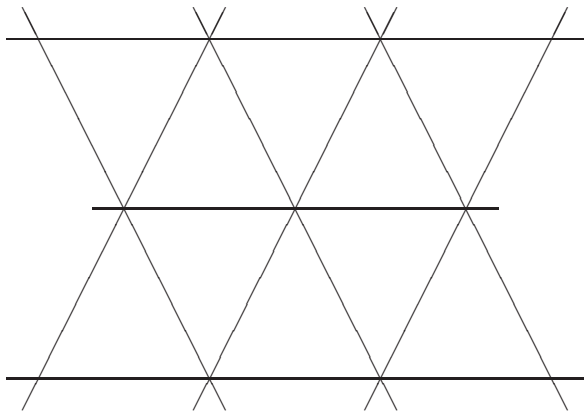
\be
\textrm{\BasicTri}
\ee
Suppose we \emph{don't} want phase to accumulate in the horizontal directions. Then, no problem, we just alternate the directions on the arrows, alternating between left pointing and right pointing (see below). This cancels out, so over reasonable distances, no phase is accumulated.
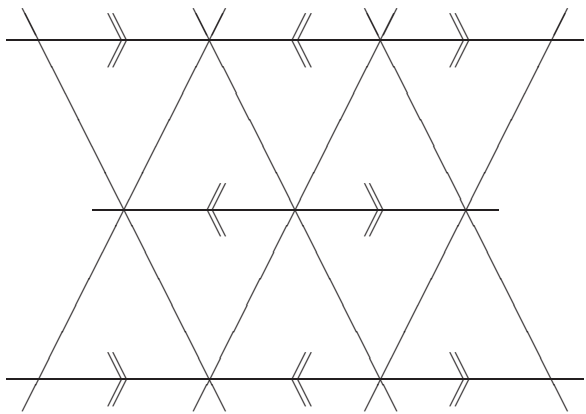
\be
\textrm{\TriAcross}
\ee
But now assume that we \emph{do} want phase to accumulate in the vertical direction. Then we simply align all the upward edges in the same direction:
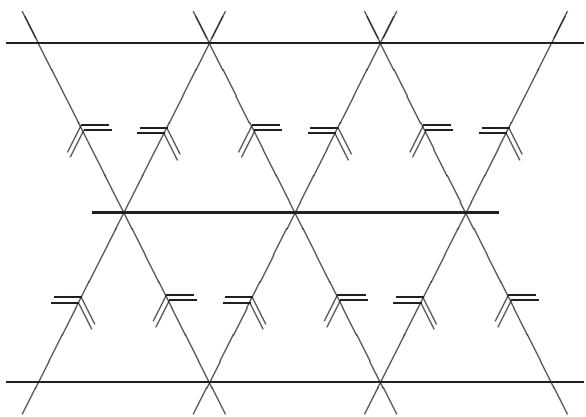
\be
\textrm{\TriUp}
\ee
And, combining the two:
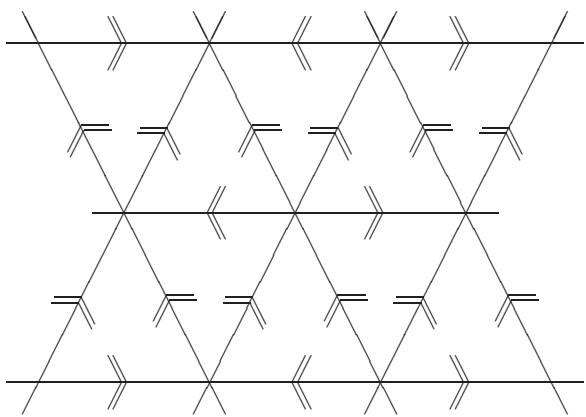
\be
\textrm{\TriAll}
\ee

We can vary the rate at which phase is accumulated in a direction, by mixing in a certain proportion of contrary arrows. And, of course, we can piece together different pieces that look locally like this to construct more complicated large scale holonomy -- letting phase accumulate along a very large loop, for instance.

\section{Higher dimensional holonomy groups}

If we have a holonomy group $G$ of dimension $d$ (and not just the one dimensional $U(1)$, as above), then we can also model this; to do so, we pick a collection $\eta_i$ of $d$ elements in the Lie algebra $\mf{g}$ of $G$, forming a basis. Then if these elements are sufficiently close to $0$, then the $\exp(\eta_i)$ will correspond to non-observable gauge changes. Furthermore each elements $g$ of $G$ will be close to an element $g'(\alpha_1,\ldots, \alpha_d)$ of type
\be
g'(\alpha_1,\ldots, \alpha_d) = \prod_{i=1}^{i=d} \left(\prod_{j = 1}^{j = \alpha_i} \exp(\eta_i)\right).
\ee
If $\epsilon_i$ are small enough, there will be an $n$ such that the $g'(\alpha_1,\ldots, \alpha_d)$ with $\sum_i |\alpha_i| \leq n$ form a subset of $G$ with a tight enough mesh, meaning that any $g$ can be well approximated by there $g'$'s.

Then we use the construction as before: with $d$ types of directed edges, each corresponding to a multiplication by $\exp(\eta_i)$, and the various `phases' (multiplications by $\exp(\eta_i)$) accumulating along paths.

If there are positive and negative distances on the network, then we shall need to use $2d$ types of edges, to distinguish between phases accumulated along positive distances and negative distances. Equivalently, have $d$ types of directed edged, and two types of undirected edges, but allow vertexes to be connected to their neighbours by two edges, not one: one from the set of undirected `distance edges' to give the sign of distance, and one from the set of directed edges to give the phase changes.

\bibliographystyle{amsalpha}
\bibliography{ref}

\end{document}